\documentclass[12 pt]{article}
\usepackage{amsmath}
\usepackage{amsthm}
\usepackage{mathrsfs}
\usepackage{amsfonts}

\usepackage[x11names]{xcolor}

\usepackage[lmargin=0.9in, rmargin=0.9in, tmargin=0.9in, bmargin=0.9in]{geometry}

\title{Complete convergence theorem \\ for a two-level contact process}
\author{Ruibo Ma \\ Peking University}

\newtheorem{tintro}{Theorem}
\newtheorem{theorem}{Theorem}[section]
\newtheorem{pro}[theorem]{Proposition}
\newtheorem{lemma}[theorem]{Lemma}

\begin {document}
\maketitle

\begin{abstract}

We study a two-level contact process. We think of fleas living on a species of animals. The animals are a supercritical contact process in $\mathbb{Z}^d$. The contact process acts as the random environment for the fleas. The fleas do not affect the animals, give birth at rate $\mu$ when they are living on a host animal, and  die at rate $\delta$ when they do not have a host animal. The main result is that if the contact process is supercritical and the fleas survive then the complete convergence theorem holds. This is done using a block construction so as a corollary we conclude that the fleas die out at their critical value. 

\end{abstract}

\section{Introduction}

While formulating our model, we think about  a population of animals with fleas that do not harm the animals while living on them. There are four states for each site in the lattice $\mathbb{Z}^d$: 0, 1, 2 and 3. The ``0'' state means the site is empty, and ``1'' means it is occupied by an animal with no fleas; ``2'' means the site is occupied by fleas without an animal, and $3=1+2$ indicates an animal with fleas. At time $t$, the state of our process is $\zeta_t: \mathbb{Z}^d  \to \{0, 1, 2, 3\}$. We also use notations $A_t = \{ x\in \mathbb{Z} ^ d : \zeta_t (x)=1\textnormal{ or }\zeta_t (x)=3 \}$ (sites occupied by animals) and $B_t = \{ x\in \mathbb{Z} ^ d : \zeta_t (x)=2\textnormal{ or }\zeta _t(x)=3 \}$ (sites occupied by fleas). Let $n_i(t,x)$ be the number of nearest neighbors of $x$ in state $i$ at time $t$. The transition rates in the model are:
\begin{align*}
0\to1\textnormal{ at rate } &   \lambda(n_1+n_3),  &  2\to3\textnormal{ at rate }  & \lambda(n_1+n_3),\\
1\to0\textnormal{ at rate } &  1, & 3\to2\textnormal{ at rate } &  1,\\
1\to3\textnormal{ at rate } &  \mu  n_3, & 2\to0\textnormal{ at rate }  & \delta.
\end{align*}
The first row corresponds to the birth of animals and the second row to the death of animals. The third row gives the birth and death rates of fleas. Notice that only fleas with an animal (a site in state 3) can give birth to a nearby site in state 1, and that fleas die only when they do not have a host (state 2). The reason for this setting is that the fleas need the resource and energy from the host animal to reproduce, and that the number of fleas on one host animal is usually large enough so that they would survive as long as the host animal is alive.

It can be seen that the rates of $0 \to 2$ and $3 \to 1$ are set as zero. They are a reflection of the sense that fleas would place their offsprings onto an available host nearby, and that fleas could live and reproduce on the body of their host animal without any risk of dying out. A natural generalization of our model could be loosening these restrictions. Further study is required for such generalization, for it drastically changes the construction of the dual (Section \ref{GR}), and we used the fact that fleas do not die out on a living host in the proof of Proposition \ref{findmu}.

A number of similar systems have been studied previously. Lanchier and Neuhauser \cite{LN1,LN2} studied a stochastic model with hosts that can be infected by a ``symbiont'' which may be a parasite that decreases the reproduction rate of the host, or a mutualist that increases it. In general, there may be several species of hosts and several species of symbionts in the system and one is interested in conditions that allow the symbionts to survive or for competing species to coexist. Motivated by ecology, they considered the situations where the symbionts  are generalists (infect all hosts) or specialists (only infect one species) and they proved a number of results about the phase diagram of the model. Later, Durrett and Lanchier \cite{DLanchier} studied the case  of a species of host with a specialist symbiont competing with another species of  host that cannot be infected. 

A few years later, Lanchier and Zhang \cite{LZ},  motivated by simulations and numerical results of Court et al.\@ \cite{Court}, studied the stacked contact process in which each site can be in state 0 (vacant), 1 (occupied by an uninfected host), or 2  (occupied by an infected host). The transition rates in the model are the following:
\begin{align*}
0\to1\textnormal{ at rate } &   \lambda_1 n_1,  &  1\to 0\textnormal{ at rate }  & 1, \\
0\to2\textnormal{ at rate } &  \lambda_1 n_2, & 2\to0\textnormal{ at rate } &  1,\\
1\to2\textnormal{ at rate } &  \lambda_2 n_2, & 2\to1\textnormal{ at rate }  & \delta.
\end{align*}
In words, the first two rows say that birth rates of individuals are not affected by their infection status, while the third row gives the rate at which infection is transferred to uninfected neighbors, and at which individuals lose their infection. Note that infected individuals, in contrast to our model,  give birth to infected ones. 

The papers mentioned above focus on showing the existence of phase transition and giving bounds on critical values. In our paper, the ``survival'' of fleas means that the probability that they never die out is positive. Here, we will show that the fleas can survive if they infect neighbors at a large enough rate (Theorem 1), but the focus of this paper will be proving the complete convergence theorem thereafter. Our methods are inspired by a paper of Remenik \cite{Remenik08}, who considered a model introduced earlier by Broman \cite{Broman}. In this system there are three states,  ``1'' (occupied), ``0'' (vacant but inhabitable), and ``$-1$''  (uninhabitable or ``blocked''). In his model, occupied sites and vacant sites become blocked at rate $\alpha$, and blocked sites become empty inhabitable sites at rate $\alpha \delta$. The unblocked sites follow a rule similar to the contact process. The phase transition was studied in the paper, and sufficient conditions for survival and extinction were presented. The paper then used a block construction similar to the one described in Section I.2 of Liggett's book  \cite{L99} to show the complete convergence theorem for the model.

At the same time as Remenik did his work, Steif and Warfheimer \cite{StWa} studied the contact process in a varying environment in which each site has an environment that can be good (1) or bad (0) and changes state according to a two-state Markov chain independent across all sites and independent of the state of the contact process. The death rate of individuals in the contact process $\delta_i$ depends  on the state of the environment but in either environment individuals give birth onto vacant sites at a constant rate. They also use a block construction similar to the one described in Liggett's book, but they use it to conclude that the contact process in their varying environment dies out at the critical value. 

In our system the animals serve as an environment where the fleas try to survive. The dynamics of our environment are given by the contact process rather than a two-state Markov chain, where the states of all sites are mutually independent, so the situation is more complicated than the ones studied in \cite{Remenik08} and \cite{StWa}. We elect to start from a configuration where the animals are distributed as an invariant measure of the contact process, so that the environment is stable for fleas. As for the fleas, we want them to occupy a finite number of sites at time $0$.

Section 2 provides basic construction of our model: the graphical representation and the dual process. While this section does not present any main results, the terminology here is crucial to the latter parts of this paper.

We describe the block construction in Section \ref{bc}. Our model has some features that can be found in the contact process, and we work from scratch and prove that the block construction commonly used in the contact process works for our model (see \cite{stflour} and \cite{L99}). Moreover,  we will describe a comparison between our model and a finite-dependent two-dimensional oriented percolation.

Recall that the two-dimensional oriented percolation was defined in \cite{stflour} using sites in the set $\mathscr  L = \{ (i, j) \in \mathbb Z^2 \mid      i+ j$  is even and $j \geq 0 \}$. From each site $(i, j)$, there are two oriented edges to $(i - 1, j + 1)$ and to $(i + 1, j + 1)$ respectively. We assign a Bernoulli random variable to each site, which indicates whether the site is open. We consider the sites that can be reached via oriented edges and open sites from the origin. The oriented percolation is called ``$m$-dependent with density at least $1-\gamma$'' if and only if for all collections $(i_k, j_k)$, $k = 1, \dotsc, n$ with $\| (i_k, j_k) - (i_l, j_l) \|_\infty > m$ as long as $k \neq l$, we have 
\begin{equation}
	P((i_k, j_k) \textnormal{ are closed for all } k) \leq \gamma ^ n  
				.
\end{equation}

Now, we briefly explain the strategy with which we prove that there is a phase transition. The details are available in Section \ref{bc}. We learn from \cite{L99} that there are properly-sized  finite space-time blocks for the contact process so that we may use them to compare the contact process to an $m$-dependent two-dimensional oriented percolation. That block construction is useful for our model as well. The block construction focuses on ``active paths'' in the graphical representation. We want these active paths to connect  the space-time blocks in a certain way. If the animals survive, we know that active paths of the animals make all those connections with high probability.  With that in mind, our idea is that if the fleas spread quickly enough to every new born animal along the active paths, then we see that the active paths of the fleas also make the same connections with high probability. Since the size of the block and the number of active paths are finite, if the birth rate of fleas is sufficiently large, we can show the forementioned result. Let $\mu_c = \inf\{ \mu\geq 0 \mid P(\textnormal{survival}) > 0 \textnormal{ for our model with parameters }\lambda, \mu, \delta \}$ which depends on $\lambda$, $\delta$ be the critical value for the parameter $\mu$. Then, we have the following theorem.

\begin{tintro}
\label{fleaperc}
 If the contact process of animals survives, then the critical value for the survival of fleas $\mu _c<   \infty$. If $\mu = \mu_c$, the fleas die out with probability 1.

\end{tintro}

Let $\nu_0$ be the points mass on the state where every site is in state $0$ and $\nu _ 1$ be the  upper invariant measure for the  animals with no fleas. Recall that the upper invariant measure for the contact process of the animals is the limiting distribution with the initial configuration of all 1's. By the complete convergence theorem of the contact process (Theorem I.2.27 in \cite{L99}), all invariant measures for the contact process are linear combinations of the upper invariant measure and the point mass on all 0's. We will define $\nu_2$, the upper invariant measure of our process with animals and fleas. We introduce the following partial order on the set $\{0, 1, 2, 3\}^{ \mathbb{Z} ^ d }$. Suppose we have two configurations $\zeta ^ 1$ and $\zeta ^ 2: \mathbb{Z} ^ d \to \{0, 1, 2, 3\}$. Define
\begin{equation}
     A^i=\{ x\in \mathbb{Z} ^ d | \zeta^i (x)=1\textnormal{ or }\zeta^i (x)=3 \}, 
\end{equation}
\begin{equation}
   B^i=\{ x\in \mathbb{Z} ^ d | \zeta^i (x)=2\textnormal{ or }\zeta^i (x)=3 \}, 
\end{equation}
for $i = 1,2$. We say $\zeta^1\leq\zeta^2$ when $A ^ 1 \subset A ^ 2$ and $B ^ 1 \subset B ^ 2$. Our process has the monotonicity property: Given two deterministic initial configurations $\zeta^1 _0\leq\zeta^2 _0$, we can construct coupled processes $\zeta^1_t$ and  $\zeta^2_t$ such that $\zeta^1 _t\leq\zeta^2 _t$, for all $ t\geq 0 $.

For probability measures on $ \{0, 1, 2, 3\}^{ \mathbb{Z} ^ d } $, we consider the ordering $\mu\leq\nu$ if and only if $\int f d\mu\leq\int f d\nu$ for all continuous increasing function $f:\{0, 1, 2, 3\}^{ \mathbb{Z} ^ d }\to\mathbb{R}$, where ``increasing'' refers to the partial order we just introduced. Again, our model has the monotonicity property that the ordering of the initial configuration is preserved through time. Suppose that $\zeta_t$ and $\xi_t$ are two realizations of our model, with $\zeta_t \sim \mu_t$, $\xi_t\sim \nu_t\forall t\geq 0$, and $\mu _ 0\leq \nu_0$. It is guaranteed that $\mu_t \leq \nu_t  \forall t\geq 0$.

We can now give the definition of the upper invariant measure. Consider $\mu _ t$ we mentioned above, and let $\mu _ 0$ be the point mass on the state of all 3's. Then $\mu _ s \leq\mu _ 0$ when $s\geq 0$, since $\mu_0$ is the largest measure on $ \{0, 1, 2, 3\}^{ \mathbb{Z} ^ d } $. Thus, $\mu _{t + s}\leq\mu _ t$ for $t, s>0$, or $\mu _ t$ is stochastically decreasing in $t$. It follows from the compactness of the set of probability measures on $ \{0, 1, 2, 3\}^{ \mathbb{Z} ^ d } $ that the limiting distribution 
\begin{equation}
\nu _2 = \lim _{t\to + \infty} \mu_t
\end{equation}
 exists. This is the upper invariant measure. We can see immediately that it is the largest invariant measure of our model.

The next section of this paper is devoted to the complete convergence theorem. We follow the convention of using superscripts to indicate the initial configuration of the process. The proof extends time to $- \infty$ and involves the dual process, which we define in Section \ref{GR}. We state the theorem below and the reader can find the detailed proof in Section \ref{cct}.


\begin{tintro}
\label{cc}
Let $\nu$, the initial configuration of our process $\zeta_t^\nu$, be such that the animals are distributed as in an invariant measure for their contact process, and the fleas occupy a deterministic finite set $B \subset\mathbb Z^d$. Let $T_1$ be the extinction time of animals, and let $T_2$ be the extinction time of fleas. We have that as $t \to \infty$,
\begin{equation}
  \zeta ^ \nu _t \Rightarrow P(T_1 < \infty) \nu_0 + P(T_1 = \infty, T_2 < \infty) \nu_1 + P(T_2 = \infty) \nu_2.
\end{equation}
\end{tintro}

\section{Graphical representation} \label{GR}

The graphical representation is a useful tool for   studying particle systems. We will now describe our graphical representation and use it to construct a dual process.

We start by introducing the graphical representation of the contact process $A_t$. We create a space-time set $\mathbb{Z} ^d \times [0, +\infty)$.     We use Poisson processes $\alpha_{x, y}(t)$ and $\beta_x (t)$ to indicate an attempted birth event from $x$ to $y$ and an attempted death event at $x$ in the contact process, respectively. The birth event processes have rate $ \lambda $ and the death event processes have rate 1, and we let them be mutually independent. Whenever an arrival of $\alpha_{x, y}(t)$ occurs, we put an arrow from $(x,t)$ to $(y,t)$ to mark a possible birth event. Whenever an arrival of $\beta_x (t)$ occurs, we put a ``$D$'' symbol at $(x,t)$ in the graphical representation, representing a possible death event. An ``active path'', as defined in \cite{L99}, goes up in time and may go through arrows in their direction, but may not cross an ``$D$'' symbol, and thus we naturally have 
\begin{equation*}
	A_t = \{ x\in \mathbb{Z} ^ d |\exists y \in A_0   \textnormal{ such that there is an active path from }(y,0)\textnormal{ to }(x,t)\}	\textnormal{.}
\end{equation*}

To define the dual process $A_s^T$, we specify a time $T > 0$ and an initial condition $A_0^T$ and go down from time $T$ searching in the graphical representation. Let $A_s^T = \{ x\in \mathbb{Z} ^ d |\exists y \in A_0^T   \textnormal{ such that there is an active path from }(x,T-s)\textnormal{ to }(y,T)\}$ so that it has the dual property
\begin{equation}
  A_T (\omega) \cap A^T_0    \neq \emptyset \Leftrightarrow A_T^T (\omega) \cap A_0\neq \emptyset .
\end{equation}
 To find $A_s^T$ in general, we start from $A_0^T \times\{ T \}$. For all $s \in  [0,T]$, the position of $A^T_s$ in the graphical representation is actually $A^T_s \times \{T - s \}$. In other words, as $s$ increases up to $T$, we move in the reverse order of time in the graphical representation.  Let the ``timeline'' of any $x \in \mathbb Z^d$ up to time $T$ be $x\times [0, T]$.   For each site in $A_s^T, 0\leq s < T$, we go down its timeline. If we encounter a ``$D$'' symbol at $(x, t-s)$, we remove the site $x$ from the dual process at time $s$.   If we see a birth arrow from $x$ to $y$ at time $t-s$, and $y \in A^T_s$,    the site $x$ is added into the dual process at time $ s $.

The sites occupied by the animals ($A_t$) for all time are given by the contact process. We choose to start the contact process from an invariant measure, so that it serves as a stable ``environment'' where the fleas grow.  The births and deaths of fleas depend on the current positions of the animals, so we may now describe the graphical representation for the fleas on top of the space-time ``map'' of the animals.

 Recall that  $B_t = \{ x\in \mathbb{Z} ^ d | \zeta_t (x)=2\textnormal{ or }\zeta_t (x)=3 \}$. Using the method just discussed, we take a realization of  $A _ t$  for all $t\geq 0$. We define some Poisson processes for births and deaths of fleas. For all ordered pairs of neighbors $(x,y)$, let $\gamma_{x, y}(t)$ be a Poisson process  with rate $\mu$. For all $x\in \mathbb{Z} ^ d$, let $\delta_x (t)$ be a Poisson process with rate $ \delta $. Naturally let all the Poisson processes be mutually independent, and independent from $A_t$. The process $\gamma_{x, y}(t)$ corresponds to possible births from $x$ to $y$, while the process $\delta_x (t)$ corresponds to possible deaths at $x$. At each arrival of $\gamma_{x, y}(t)$, the fleas at $x$ give birth to the site $y$ if animals are present at both $x$ and $   y$. At each arrival of $\delta_x (t)$, the fleas at $x$ die, if there are fleas and no animal at site $x$.

All these events are marked by arrows and symbols in $\mathbb{Z} ^d \times [0, +\infty)$. Whenever an arrival of $\gamma_{x, y}(t)$ occurs, we draw an arrow from $(x,t)$ to $(y,t)$ in the graphical representation if $x, y\in A_t$, representing the possible birth of the fleas at $y$. Whenever an arrival of $\delta_x (t)$ occurs, we put an ``*'' at $(x,t)$ in the graphical representation if $x$ is not occupied by an animal at time $t$, representing the death of fleas, if there are any.

We now can determine the state of any site at any time in the graphical representation.   An ``active path'' for fleas goes up in time and may go through  arrows of fleas in the correct direction, but is not allowed to cross an ``*''. A site $x$ is occupied by fleas at time $t$ if and only if there is an active path from $(y,0)$ to $(x,t)$ where $y$ is occupied by fleas in the initial configuration. This leads to the definition of the dual process $B^T_s$. We first specify a $D\subset\mathbb{Z}^d$ as $B^T_0$. Let $B_s^T = \{ x\in \mathbb{Z} ^ d |\exists y \in D \textnormal{ such that there is an active path from }(x,T-s)\textnormal{ to }(y,T)\}$, $0\leq s\leq T$. To visualize how to find $B^T_s$, we start from $B^T_0\times\{T\}$. We are going back in time along the timeline of each member $x$ of $B^T_0$. When we encounter an arrow from $(y,T-s)$ to $(x, T-s)$, $y$ is added to $B^T_s$ at time $s$, and we also go down the timeline of the new member $y$. When we encounter an asterisk at $(x,T-s)$, $x$ is removed from $B^T_s$ at time $s$, and we disregard the timeline of $x$. $B^T_s$ has the desired dual property that with the whole history of animals given, then $B_T\cap D\neq \emptyset$ if and only if $B^T_T\cap B_0\neq\emptyset$. We formalize that property in the following proposition.




\begin{pro}\label{dual}
Let $T > 0$, and let $B, D \subset \mathbb Z ^d$ be  finite sets. If $B_0 = B$, and $B^T_t$ is a dual process of fleas with the initial configuration $B^T_0 = D$, we know that
\begin{equation}
         \{\omega : B_ T (\omega)  \cap D \neq \emptyset \} =  \{\omega : B^T_T (\omega)  \cap B\neq \emptyset\}
 .
\end{equation}
\end{pro}
\begin{proof}
For each $\omega \in \Omega$, we have information from those Poisson processes about the behavior of animals and fleas. Thus, we can determine the ``active paths'' for fleas accordingly. By the definition of the dual process, $B^T_T\cap B\neq \emptyset$ iff there exist $x \in B, y \in D$ such that there is an ``active path'' from $(x, 0)$ to $(y, T)$ in the graphical representation. Since the fleas give birth along the active paths, that is exactly $B_ T    \cap D \neq \emptyset$.
\end{proof}

\section{Block construction and proof of Theorem 1}
\label{bc}
This part of the paper provides   a block construction argument. The techniques we utilize have been present in the field for long, and they are available in the 1990 paper \cite{BG1990} and Section I.2 of Liggett's 1999 book \cite{L99}. The main result of this section, Theorem \ref{block_5}, gives equivalent conditions of survival of the process.

Recall that  superscripts seen next to stochastic processes  indicate the initial configuration, as we stated before Theorem \ref{cc}. From now on, to simplify symbols involving the fleas, we may indicate only the sites initially occupied by fleas, and we omit the invariant measure regarding the animals. This is acceptable, since the measure of the animals does not influence the following proofs. For instance, we may use $B^B_t$ to denote the set occupied by fleas at time $t$ with the condition that the fleas occupy $B \subset \mathbb Z^d$ at time $0$, and we may use $B^{ T, D}_t$ for a dual process of fleas, with the condition that $B^T_0 = D$, etc.

\begin{lemma}\label{block_1}
Suppose the fleas survive. Then 
\begin{equation}
  \lim_{n\to\infty}P  (B ^ { [-n, n] ^ d } _t\neq\emptyset\,\,\forall t\geq0) =1.     
\end{equation}
  
\end{lemma}

\begin{proof}
For a positive integer $k$, we let $X _ k(\omega)=1_{\{ B ^ {\{ (k,0,...,0) \}}_t\textnormal{ survives}\}}(\omega)$. $\{X_k \}$ is then a stationary sequence. By Birkhoff's ergodic theorem, $\frac 1 n(X _ 1+ ... +X _ n) \to E(X_1 | I)$ a.s., where $I$ is the invariant $ \sigma $-field. Note that an invariant event with respect to $\{ X _ k \}$ is also an invariant event with respect to all the Poisson processes $\alpha_{x, y},\beta_x,\gamma_{x, y},\delta_x$ under a shift. Since the Poisson processes are all i.i.d., $I$ is trivial by L\' evy's zero-one law.

Now, $\frac 1 n(X _ 1+ ... +X _ n) \to E(X_1) >0$ a.s. This implies $M=\inf \{k|X_k = 1\} < \infty$ a.s. For any $\varepsilon > 0$, there exists a constant $N > 0$ such that $P (M<N) > 1 - \varepsilon $. This implies that $P (B_t ^{ [-N, N] ^ d }      \neq\emptyset\,\,\forall t\geq0)>1-\varepsilon$, and  the proof is complete.
\end{proof}

We further define $_LA_t$, $_LB_t$ as the truncated process of animals and fleas, respectively, where births of animals and fleas only originate from a site inside $(-L, L) ^ d$, and the vertical segments of all active paths are within $(-L, L)^ d \times [0, \infty)$ in the graphical representation. The next few lemmas lead to the main result of this section.

\begin{lemma}
\label{block_2}
For any finite $A\subset\mathbb{Z}^d$ and any $N\geq1$, 
\begin{equation}
 \lim _{t\to\infty}\lim _{L\to\infty} P (|_LB^A _t|\geq N)= P  (B^A _t\neq\emptyset\,\,\forall t\geq0).
\end{equation}
\end{lemma}

\begin{proof}
By the monotonicity of the process, we see that the limit as $ L\to\infty $ is $P(|B ^ A_t|\geq N)$. We then suppose that we have fleas at $ N $ sites at time $ s $. To guarantee that the fleas die out, we only need the following set of events: The death of an animal and the death of fleas happen in that order at the $N$ sites before any spreading of animal or fleas can occur to or from these sites. Thus the conditional extinction probability of the fleas can be bounded from below by a number which only depends on $ N $. By the martingale convergence theorem, 
\begin{equation}
P (| B^A_t |=0 \textnormal{ for some }t>0 | \mathscr{F}_s ) \to 1 _ {\{| B_t^A |=0 \textnormal{ for some }t>0 \}}
\end{equation}
almost surely as $s\to +\infty$, where $| B_t^A |$ is the number of sites with fleas at time $ t $, and $\mathscr{F}_s$ is the $ \sigma $-field generated by all the Poisson processes up to time $ s $ and the initial condition. For any $\omega$ such that $\liminf \limits_{s\to \infty} |B_s^A(\omega)| < \infty$, we must see that $\limsup \limits _{s\to \infty} P (| B^A_t |=0 \textnormal{ for some }t>0 | \mathscr{F}_s )(\omega)$ is strictly greater than 0. It follows that $\omega \in\{| B_t^A |=0 \textnormal{ for some }t>0 \}$, and that $| B^ A _t |\to\infty$ a.s.\@ on $\{\omega\mid | B_t^A  |\neq   0\textnormal{ for all }t>0 \}$. Seeing that
\begin{equation}
 \liminf _{t\to\infty}\lim _{L\to\infty} P (|_LB^A _t|\geq N)\geq    P  (B^A _t\neq\emptyset\,\,\forall t\geq0)
\end{equation}
for $N\geq 1$ by the arguments above and that
\begin{equation}
     \limsup _{t\to\infty}\lim _{L\to\infty} P (|_LB^A _t|\geq N) \leq    P  (B^A _t\neq\emptyset\,\,\forall t\geq0)
   ,
\end{equation}
the result follows.
\end{proof}

Let $S(L, T) = \{x \in \mathbb{Z} ^ d: \|x\|_\infty = L\} \times[0,T]$, where $\| \cdot \|_\infty$ is the $l^\infty$-norm, and let $N(L, T)$ be the maximal number of points in a subset $F$ of $S(L,T)$ with the properties (i) $(x, s)\in   F$ implies $x\in{_LB_s}$; (ii) if $(x,s'), (x,s'')\in F$, then $|s' -s''|\geq1$. Let $S_+ (L,T) = (\{L\}\times   \{0, ..., L\}^{d-1})\times[0,T]$, and let $N_+(L, T)$ be the maximal number of points in a subset $F$ of $S_+(L,T)$ with the same properties.

The next lemma is proved by applying  positive correlation to increasing functions.

\begin{lemma}
\label{block_3}
Suppose that at time 0, the animals are in their invariant measure, and the fleas occupy $[-n, n] ^ d$. Then we have
\begin{equation}
   P ( |_LB_T\cap [0,L]^d | \leq N)\leq P ( |_LB_T|\leq 2^dN)^{2^{-d}}
\end{equation}
for $N\geq 1,L \geq n$ and
\begin{equation}
 P (N_+(L,T)\leq M)^{d2^{d}}\leq  P (N(L,T)\leq Md2^d)
\end{equation}
for  $ L >n$, $ T >0$, $ M\geq 1 $.

\end{lemma}
\begin{proof}
There are $2 ^ d$ orthants in the space $\mathbb{Z} ^ d$. Let $X _ k$ be the size of the intersection of $_LB_T$ and the $ k $th orthant. Now the $X _ k$'s are identically distributed, and the positive correlation of our model applies. (See Corollary B18 and Proposition I.2.6 in \cite{L99}.) Thus, we have
\begin{align}
     \nonumber
P ( |_LB_T|\leq 2^dN) &\geq P (X_k \leq N, 1 \leq k \leq 2^d)\\
 &\geq P ( |_LB_T\cap [0,L]^d | \leq N) ^ {2^d}.
\end{align}
This proves the first inequality. The second can be proved by similar means.
\end{proof}
\begin{lemma} \label{block_4}
Suppose that at time 0, the animals are in an invariant measure, and the fleas occupy a finite set $A \subset \mathbb{Z} ^ d$, and suppose $L_j \uparrow \infty$ and $T_j \uparrow \infty$. For any $M,N \geq 1$,
\begin{equation}
     \limsup _ {j\to\infty} P ( N(L_j,T_j)\leq M ) P ( |_{L_j} B_{T_j}|\leq N )\leq P (| B_s |=0 \textnormal{ for some }s>0 )  
      .
\end{equation}
\end{lemma}

\begin{proof}
As in Lemma \ref{block_2},      we bound the probability of extinction from below, given $A \subset (-L, L) ^ d$, and $ |_LB_T|+N(L,T)\leq k$. For all sites in $ _LB_T $, we want the death of the animal on the site before any birth of  fleas to a neighboring site, and then the death of fleas before any animal births onto this site. For different sites in $_LB_T$, these events are mutually independent.

For the $N(L, T)$ space-time points on the side of the space-time box, we want to avoid them giving birth to any neighboring sites in the nontruncated process. The total length of intervals where these sites are occupied by fleas is no more than $2N(L, T)$. So the probability that no birth of fleas occurs along these intervals can be bounded from below using $N(L, T)$. To sum up, the conditional extinction probability can be bounded from below by a number which only depends on $ k $.

Let $G=\{| B_s |=0 \textnormal{ for some }s>0\}$, and $H _j=\{ |_{L_j}B_{T_j}|+N(L_j,T_j)\leq k\}$ for a fixed $ k $. By the martingale convergence theorem,
\begin{equation}
   P ( G | \mathscr{ F }_{L _ j,T _ j} ) \to 1 _ G \textnormal{ a.s.}
\end{equation}
as $j\to\infty$, where $ \mathscr{ F }_{L _ j,T _ j} $ is the $\sigma$-field generated by all the Poisson processes involved in $_{L_j}B_t$ up to time $T_j$. Now $P ( G | \mathscr{ F }_{L _ j,T _ j} )$ is bounded from below on $H _j$. Therefore, $\{H _j\textnormal{ i.o.}\}\backslash G$ has measure 0. Thus, $\limsup _{j\to\infty}P  (H _j) = \lim_{n\to\infty} \sup_{j \geq n} P(H_j)   \leq \lim_{n\to \infty} P\left(\bigcup_{j\geq n} H_j \right) =P(H _j\textnormal{ i.o.})\leq P(G) $. Based on positive correlations, for $L\geq 1, T>0, M, N \geq 1$,   
\begin{align*}
P ( |_LB_T|+N(L,T)\leq M + N)&\geq P   (N(L,T)\leq M,|_LB_T|\leq N)   \\
   &\geq P(N(L,T)\leq M) P ( |_LB_T|\leq N)   .
\end{align*}
Thus, after we let $k = M+N$, we have
\begin{align*}
\limsup _ {j\to\infty} P ( N(L_j,T_j)\leq M ) P ( |_{L_j} B_{T_j}|\leq N ) & \leq  \limsup _ {j\to\infty}   P (   |_{L_j} B_{T_j}| + N(L_j,T_j)\leq M+N ) \\
  &= \limsup _{j\to\infty}P  (H _j)\\
  & \leq P(G) \\
  & = P (| B_s |=0 \textnormal{ for some }s>0 )   .
\end{align*}
\end{proof}

With these lemmas, we now prove the main theorem of this section.

\begin{theorem}
\label{block_5}
The fleas survive if and only if for any given $\varepsilon > 0$, there are $n, L, T >0$ which satisfy the following conditions:
\begin{equation}
 P  (_{2n+L}B ^ { [-n, n] ^ d }_ {T+1}\supset (x+[-n,n]^{d})\cap \mathbb{Z} ^ d \textnormal{ for some }x\in[0,L)^d) \geq1-\varepsilon
 \label{Ptop}
\end{equation}
and
\begin{equation}
     \begin{aligned}
&     P (_{2n+L}B ^ { [-n, n] ^ d }_ {t+1}\supset (x+[-n,n]^{d})\cap \mathbb{Z} ^ d \textnormal{ for some }0\leq t\leq T     \textnormal{ and some }
    \\      
 &\qquad\qquad\qquad\qquad\qquad\qquad\qquad\qquad           x\in\{L+n\}\times [0,L)^{d - 1})\geq1-\varepsilon.
     \end{aligned}
\label{Pmove}
\end{equation}
\end{theorem}

\begin{proof}
The necessity can be proved using Lemmas \ref{block_1}-\ref{block_4}. Rather than copying the classic proof of Theorem I.2.12 of Liggett's book \cite{L99}, we refer the reader to that book for a detailed proof.

To show  the sufficiency, we can assume that the process starts with $[-n, n] ^ d$ occupied by fleas. That is fine since there is a positive probability that  $B_1^{ \{ 0 \}} \supset [-n, n] ^ d$. By Proposition I.2.20, Proposition I.2.22 and Theorem 2.23 of \cite{L99}, the fleas survive.
\end{proof}

The next proposition leads to Theorem \ref{fleaperc}.

\begin{pro}
\label{findmu}
Let $\lambda>0$ be supercritical and let $\delta>0$ be fixed. For any $\varepsilon>0$, we can find a $\mu < \infty$ such that the conditions in Theorem \ref{block_5} are satisfied with some $n, L, T$.
\end{pro}

\begin{proof}
In the supercritical case, the space-time conditions of the animals are satisfied, namely
\begin{equation}
P  (_{2n+L}A _ {T+1}^ { [-n, n] ^ d }    \supset (x _1+[-n,n]^{d})\cap \mathbb{Z} ^ d \textnormal{ for some }x _1\in[0,L)^d) \geq1-\varepsilon ,
\end{equation}
and
\begin{equation}
   \begin{aligned}
P ( & _{2n+L}A _ {t+1}^ { [-n, n] ^ d } \supset (x _2+[-n,n]^{d})\cap \mathbb{Z} ^ d \textnormal{ for some }0\leq t\leq T     \textnormal{ and some }          \\
 &\qquad\qquad\qquad\qquad\qquad\qquad\qquad\qquad
  x _2\in\{L+n\}\times [0,L)^{d - 1})\geq1-\varepsilon.
  \end{aligned}
\end{equation}

In the graphical representation, the conditions translate to active paths from $([-n,n] ^ d\cap\mathbb{Z}^d)\times\{ 0 \}$ to $( (x_1+[-n,n]^{d})\cap \mathbb{Z} ^ d )\times\{ T + 1 \}$ and $( (x_2+[-n,n]^{d})\cap \mathbb{Z} ^ d )\times\{ t + 1 \}$. There is high probability of finding these paths to all sites in the moved cube ($(x_k+[-n,n]^{d})\cap \mathbb{Z} ^ d$, $k = 1, 2$) at the corresponding time ($t + 1$ or $T + 1$).  Let $\Gamma_1$ be the collection of paths from $([-n,n] ^ d\cap\mathbb{Z}^d)\times\{ 0 \}$ to $( (x_1+[-n,n]^{d})\cap \mathbb{Z} ^ d )\times\{ T + 1 \}$, and let $\Gamma_2$ be the collection of the paths from $([-n,n] ^ d\cap\mathbb{Z}^d)\times\{ 0 \}$ to $( (x_2+[-n,n]^{d})\cap \mathbb{Z} ^ d )\times\{ t + 1 \}$. Even when either event does not occur, $\Gamma _ 1$ and $\Gamma _ 2$ are still well defined.

We want to bound the number of births on all these paths, so that we know the number of births the fleas have to give to spread themselves through the same paths. The total number of jumps in $\Gamma _ 1$ and $\Gamma _ 2$ are finite with probability 1. Take an $N>0$ so that the probability that there are more than $N$ jumps in either $\Gamma _ 1$ or $\Gamma _ 2$ is smaller than $\varepsilon$.

We also want to make sure there is enough time for the fleas to give birth. After each birth on the animal paths, if the fleas also give birth onto the new born animal before the end of the path, the next birth on the path and the death of either animal (parent or offspring), they are able to spread themselves on the same path. We call the length for them to give birth  the ``birth window''. The smallest birth window among all paths in $\Gamma _ 1$ and $\Gamma _ 2$ is positive with probability 1. Let $\omega > 0$ be such that the probability that the smallest window among all paths in $\Gamma _ 1$ and $\Gamma _ 2$ is less than $\omega$ is smaller than $\varepsilon$.

Now with high probability, we have two collections of desired paths of host animals such that we see no more than $N$ birth events in either one, and each ``birth window''  is at least $\omega$. By independence, on that event, we can find a finite $ \mu $ so that the fleas give birth during each birth window with probability at least $1-3\varepsilon$. Thus, the probabilities in \eqref{Ptop} and \eqref{Pmove} are at least $1-4 \varepsilon$. Since $\varepsilon$ is arbitrary, the conditions for fleas are proved.
\end{proof}

The seemingly unpleasant flaw of this proposition is that $ \mu $ depends on the choice of $\varepsilon$. It may seem at first that it is not enough to show the survival of fleas for a fixed $\mu$. However, knowledge of the oriented percolation eliminates the need to worry.

\begin{proof}[Proof of Theorem \ref{fleaperc}]
When the animals are supercritical, Proposition \ref{findmu} shows dominance over an $m$-dependent oriented percolation under certain condition. To show survival of fleas, we need the $m$-dependent oriented percolation to be supercritical. Hence, we  let $\varepsilon$ be small such that all $m$-dependent oriented percolations with density at least $1 - \varepsilon$ survive. By Proposition \ref{findmu},  there exists $\mu_c<\infty$ such that whenever  $\mu>\mu_c$, the fleas survive.

Since the conditions \eqref{Ptop} and \eqref{Pmove} are defined in a finite space-time set, we see that $\mu - \varepsilon_0$ is supercritical for all sufficiently small $\varepsilon_0 > 0$, when $\mu$ is supercritical. With $\lambda, \delta$ fixed, the set of all supercritical values of $\mu$ is an open  set. Hence, the value $\mu_c$ must be subcritical, which proves the rest of this theorem.
\end{proof}

\section{Complete convergence}
\label{cct}

We are ready to use the block construction to prove the complete convergence theorem. By the complete convergence theorem of the contact process, we reduce the desired result, Theorem \ref{cc}, to the complete convergence of fleas, i.e.\@ $B^\nu _t\Rightarrow P (T_2 < \infty) \underline{\nu} + P  (T_2 = \infty) \bar{\nu}$, where $ \underline{\nu} $ and $ \bar{\nu} $ are the distributions of fleas under measures $\nu_0$ and $\nu_2$, respectively.

We first need to describe the block construction for the dual process.

\subsection{Block construction for the dual process}

We need to rebuild the graphical representation to provide the block construction result for the dual process. Since the dual goes back in time, we need to extend the time of the graphical representation to $-\infty$ so that the time of the dual process extends to $\infty$. That is done by the Kolmogorov  extension theorem.

For all $n\geq 1$ and $t_1, \dotsc, t_n \in  \mathbb R$, we consider a contact process of animals starting from time $\min \{t_1, \dotsc, t_n \}$. As we did before, the initial configuration is an invariant measure of animals. Let the distribution of $(A_{t_1}, \dotsc, A_{t_n})$ be $\nu_{t_1, \dotsc, t_n}$. It is natural that all such measures are consistent, and the Kolmogorov extension theorem applies. Therefore, the existence of the contact process on $(- \infty, \infty)$ is established. We may extend the Poisson processes we defined in Section \ref{GR} to time $-\infty$ so that the symbols representing possible birth and death events are properly marked in the negative-time section of the graphical representation. Now that the behavior of animals on $(-\infty, \infty)$ is known, and all symbols of flea events are present, the dual process of fleas $B^T_s$ can be defined for all $s\geq 0$.

Since the dual process runs up to time $\infty$, we can see that Lemmas \ref{block_1}-\ref{block_4} apply.  Lemmas \ref{block_1} and \ref{block_3} hold  for the dual process by the same proof. For the other two lemmas, the key step of the proof was to bound the probability of dying out from below using the size of the occupied set. That can be done for the dual process as well. Suppose that $B^T_s$ has size $N$ for some $s$. We see that the dual process dies out if for each of the $N$ sites, we encounter an animal birth event (note that we are going back in time) and a flea death symbol before any other symbols. Thus, the extinction probability  can be bounded from below by a positive constant only depending on $N$. 

That argument means that the proofs of Lemmas \ref{block_2} and \ref{block_4} apply to the dual process as well. Hence we have the following result.

\begin{pro}
The dual process $B^{T_0}_t$ survives if and only if for any given $\varepsilon > 0$, there are $n, L, T >0$ which satisfy the following conditions:
\begin{equation}
    P (_{2n+L}B ^{T_0, [-n, n] ^ d}_ {T+1}\supset (x+[-n,n]^{d})\cap \mathbb{Z} ^ d \textnormal{ for some }x\in[0,L)^d) \geq1-\varepsilon
\end{equation}
and
\begin{equation}
     \begin{aligned}
     P  (  &_{2n+L}B ^{T_0, [-n, n] ^ d     }_ {t+1}\supset (x+[-n,n]^{d})\cap \mathbb{Z} ^ d \textnormal{ for some }0\leq t\leq T      \textnormal{ and some }     \\
      &\qquad\qquad\qquad\qquad\qquad\qquad\qquad\qquad
      x\in\{L+n\}\times [0,L)^{d - 1})\geq1-\varepsilon.
   \end{aligned}
\end{equation}
\end{pro}

\subsection{A helpful lemma}

Leading to the proof of Theorem \ref{cc}, we use the graphical representation in a slightly different way. We pick a large $T>0$, two finite sets $B, D \subset \mathbb Z^d$, and we run the process of fleas $B^B_t$ and the dual process $B^{T,D} _s$ both up to time $\frac T 2$.

The block construction for the dual process leads us to the following lemma, which is key to the proof of Theorem \ref{cc}.

\begin{lemma}\label{cc1}
Suppose that $B, D\subset \mathbb{Z}^d$ are finite and nonempty, and that $\varepsilon_0>0$. There is a $T_0>0$ such that when $T\geq T_0$, 
\begin{equation}
     P(B^B _\frac{T}{2}\neq \emptyset, B^{T, D} _\frac{T}{2}\neq \emptyset, B^B_\frac{T}{2}\cap B^ {T, D}     _\frac{T}{2}=\emptyset)\leq \varepsilon_0
     .
\label{bad_event_bd}
\end{equation}
\end{lemma}

\begin{proof}
The proof is done in a few steps. Fix an $0<\varepsilon < 1/6$. 

\textbf{Step 1}. Recall that the block construction provides a comparison with an $m$-dependent oriented percolation, defined on page 138 of \cite{stflour}. By Theorem B26 of \cite{L99}, this oriented percolation with sufficiently large density dominates the original version with any density $<1$ defined in \cite{D1984}.

Suppose that the conditions \eqref{Ptop}, \eqref{Pmove} hold (for all $\varepsilon > 0$, not to be confused with the fixed $\varepsilon$ above), and $\varepsilon_1>0$. By Proposition I.2.22 of \cite{L99}, there exists $n, a, b>0$ such that for $\forall (x, s) \in [-a, a] ^d\times [0, b]$, the probability that  there are active paths of fleas, within $[-5a, 5a]^d \times [0, 6b] $, starting from $(x, s) + [-n, n]^d \times \{0\} $ to every point in $(y, t) + [-n, n]^d \times \{0\} $ for some $(y, t) \in [a, 3a] \times [-a,a] ^ {d-1} \times [5b,6b]$ is greater than $1 - \varepsilon_1$.

Let  $\mathscr  L = \{ (i, j) \in \mathbb Z^2 \mid      i+ j$  is even and $j \geq 0 \}$. We say that  $(i, j)$  is open if $B^{[-n, n]^d}_t \supset x + [-n,n]^d$ for some $(x, t) \in [(2i-1)a, (2i+1)a] \times [-a, a]^{d-1} \times [5jb, (5j+1)b ]$. By Proposition I.2.22 of \cite{L99},  the stochastic process of open vertices dominate an  $ m $-dependent oriented percolation with density $1 - \varepsilon_1$.

We consider the $ m $-dependent oriented percolation with density $1 - \varepsilon_1$. Recall that the notations $l_n$ and $r_n$ for the independent oriented percolation are defined as
$$
l_n = \inf \{x : (x, n) \in \mathscr L , \textnormal { and it can be reached from } (0, 0) \textnormal { via open sites and oriented edges} \},
$$
$$
 r_n = \sup \{x : (x, n) \in \mathscr L ,  \textnormal { and it can be reached from } (0, 0) \textnormal { via open sites and oriented edges} \}   
,
$$
and recall that the oriented percolation ``survives'' if infinitely many vertices can be reached from the origin. By    Section 13 of \cite{D1984}, the proportion of sites in $\mathscr L$ between $(l_n, n)$ and $(r_n, n)$ that can be reached from $(0, 0)$ converges to the probability of survival, denoted by $\rho$, almost surely on the event of survival. Moreover, by the proof presented in the paper \cite{D1984}, we see that the same holds if we replace $l_n$ by $l_n '$, and $r_n$ by $r_n '$, $l_n\leq l'_n \leq r'_n\leq r_n$, and $\lim_{n\to\infty} \frac{l'_n}{n}$, $\lim_{n\to\infty} \frac{r'_n}{n}$ exist. We may pick $\varepsilon_1$ properly such that any $ m $-dependent oriented percolation with density at least $1 - \varepsilon_1$ dominates an independent oriented percolation  with survival probability $     1 - \varepsilon > 5/6$. It implies that there is an $N_1 \geq 0$ such that the probability that the $ m $-dependent oriented percolation with density  $1 - \varepsilon_1$   survives but less than $2/3$ of the sites between $(l' _n, n)$ and $(r' _n, n)$ can be reached from the origin for some $n > N_1$ is less than $\varepsilon$.

\textbf{Step 2}. Consider the block construction in Step 1. With the $\varepsilon_1$ chosen in Step 1, we obtained constants $n, a, b > 0$ for the process of fleas. With the same   $\varepsilon_1$, there are constants $\tilde n, \tilde a, \tilde b > 0$ for the block construction of the dual.

Recall that in the initial configuration of the block construction,  a cube around the origin is occupied by fleas.    Instead, the initial configuration in this lemma is general, and we can show that a cube of the same size exists (possibly centered elsewhere) with high probability given that the fleas survive sufficiently long. It is known that the probability that the fleas survive up to time $t$ but we do not see  an occupied cube of the same size anywhere by time $t$ goes to 0 as $t\to +\infty$.  Thus, there is a $t_1>0$ such that the probability that the fleas survive until time $t_1$ without forming a cube of size $2n$ is less than $\varepsilon$. For the same reason, there is a $t_2>0$ such that the probability that the $B^ {T, D} _ s$ survives until time $t_2$ without forming a cube of size $2 \tilde{n}$ is less than $\varepsilon$ for any $T>0$. Since the sets $B, D$ are finite, we can find a number $M>0$ such that the probability that the cubes within $B^B_t$, $B^ {T, D} _s$ are formed by time $t_1, t_2$, respectively and at least one of their centers is out of $[-M,M]^d$ is less than $\varepsilon$.

\textbf{Step 3}. In this step, we consider the event where the two big occupied cubes are formed and centered within $[-M, M]^d$ in time. All probabilities mentioned in this step refer to the probability of an intersection with this event.  After the two cubes are occupied, we couple the processes $B^B_t$ and $B^ {T, D} _s$ with the dominated $m$-dependent oriented percolation. 

Suppose that $x_1 + [-n, n]^d \subset B_{s_1} ^B$, with $s_1\leq t_1$,  that $x_2 + [-\tilde n, \tilde n]^d \subset B_{s_2}^{T,D}$, with $s_2 \leq t_2$, and that $x_1, x_2 \in [-M, M]^d$.      Let $M_1 = 24(a+\tilde a)$, and $M_2 = M + 5(a+ \tilde a     )$, and let $n_1 = [(T/2 - s_1) / (5b)] - 1$, $n_2 = [(T/2 - s_2) / (5     \tilde{b} )]-1$, where $[x]$ is the floor function. For the fleas, we let $(0,0)$ in the $m$-dependent oriented percolation (denoted by $\xi _n^{(1)}$) correspond to the space-time region $(x_1 + [-a, a]^{d})  \times [s_1, s_1+b]$ in the graphical representation. Moreover, the vertex $(i, j) \in \mathscr L$ of the $m$-dependent oriented percolation correlates with the region $(x _ 1 + [(2i-1)a, (2i+1)a]\times [-a, a]^{d-1}) \times [s_1 + 5jb, s_1 + (5j+1)b]$. We run the oriented percolation until time $n_1$, which means that we consider the space-time regions of $(i, j)$ with $j \leq n_1$ and $-j \leq i \leq j$, because the oriented percolation spreads out by at most distance $j$ at time $j$. The region of $(i, n_1), -n_1 \leq i \leq n_1$ is each contained in or intersects with a region $R_ k = [kM_1, (k+1) M_1] \times [-M_2, M_2] ^ {d-1} \times [s_1 + 5 n_1 b, s_1 + (5n_1+1) b   ]$ for some $k \in \mathbb{Z}$ in the graphical representation.  

Likewise, every site of the $m$-dependent oriented percolation coupled with the dual process, denoted by $\xi_n^{(2)} $, correlates with a space-time region. The root $(0, 0)$ correlates with space-time coordinates $(x_2 + [- \tilde a, \tilde a]^{d})  \times [T - s_2 - \tilde b, T - s_2] $. Generally, $(i, j) \in \mathscr L$ corresponds to $(x _ 2 + [(2i-1) \tilde a, (2i+1) \tilde a]\times [-\tilde a, \tilde a ]^{d-1}) \times [T - s_2 - (5j + 1)\tilde b, T - s_2 - 5j \tilde b]$, $   0 \leq j \leq n_2$.  When $j = n_2$, each region of $(i, n_2),   -n_2 \leq i\leq n_2$ is  contained in or intersects with a region $\tilde{R} _k=[kM_1, (k+1) M_1] \times [-M_2, M_2] ^ {d-1} \times [T - s_2 - (5 n_2+1) \tilde b, T - s_2 - 5n_2 \tilde b]$ for a $k \in \mathbb{Z}$ in the graphical representation.  

Let $l_n$, $r_n$ be defined as in Step 1 for $\xi_n^{(1)}$, and   let $\tilde l_n, \tilde r_n$ be defined likewise for $\xi_n^{( 2 )}$. Let $K$ be the set that contains all values of $k$ such that $R_k$ contains a region of a site between $(l_{n_1}, n_1)$ and $(r_{n_1}, n_1)$ of $\xi_n^{(1)}$, and that $\tilde R_k$ contains a region of a site between $(\tilde  l_{n_2}, n_2)$ and $(\tilde r_{n_2}, n_2)$ of $\xi_n^{( 2 )}$, consider the vertices $(i, n_1)$ of $\xi_n^{(1)}$ such that its region is contained in $\bigcup_K R_k$, and consider the vertices $(i, n_2)$ of $\xi_n^{(2)}$ such that its region is contained in $\bigcup_K \tilde{R}_k$. As long as $T$ is sufficiently large, the probability that $\xi_n^{(1)}$ survives until time $n_1$ but less than $2/3$ of the sites within the forementioned range are in $\xi_{n _  1} ^{(1)}$ is less than $\varepsilon$, and the probability that $\xi_n^{(2)}$ survives until time $n_2$ but less than $2/3$ of the sites within the forementioned range are in $\xi_{n _2} ^{(2)}$ is less than $\varepsilon$. Thus, with high probability, the proportion of sites in the range of $\xi_n^{(1)}$ that are open with their region contained in an $R_k$, $k \in K$ is at least $7/12$, and the proportion of sites in the range of $\xi_n^{(2)}$ that are open with their region contained in an $\tilde R_k$, $k \in K$ is at least $7/12$,  since the proportion of those regions lapping over $R_k$ and $R_{k+1}$ or $\tilde{R}_k$ and $\tilde{R}_{k+1}$ is at most $1/12$. It implies that the proportion of regions $R_k$, $k \in K$ that contains a region of $(i, n_1)$, $i \in \xi_{n _  1} ^{(1)}$ is at least $7/12$ with high probability, and that  the proportion of regions $\tilde {R}_k$, $k \in K$ that contains a region of $(i, n_2)$,  $i \in \xi_{n _  2} ^{(2)}$ is at least $7/12$ with high probability. Thus, among all members $k \in K$, at least $1/6$ satisfies that $R_k$ contains a region of $(i, n_1)$, $i \in \xi_{n _  1} ^{(1)}$, and that $\tilde {R}_k$ contains a region of $(i, n_2)$,  $i \in \xi_{n _  2} ^{(2)}$ with high probability. Let $K^\prime$ be the set of such $k$.


\textbf{Step 4}. Now consider the regions $R_k$ and $\tilde R_k$ for $k \in K^\prime$. Given the condition at the end of the last step, we may bound the probability that $B^B_\frac{T}{2}$ intersects with $ B^ {T, D}     _\frac{T}{2}$ in $[kM_1, (k+1) M_1] \times [-M_2, M_2] ^ {d-1}$ from below. Consider 
\begin{align}   
\notag
&   \inf _{x, y \in [ 0,  M_1 ] \times [-M_2, M_2] ^ {d-1}, 4b + 4\tilde b \leq t \leq 10b + 10\tilde b} P(  \textnormal{there is an active path in}\\
     \label{intersect inf}
&   \quad\quad\quad\quad\quad\quad\quad\quad\quad\quad      [ - (n + \tilde n),  M_1 + (n + \tilde n)] \times [-M_2, M_2] ^ {d-1} \times[0, t]  \\
 & \quad\quad\quad\quad\quad  \textnormal{     from an } (x', 0) \in (x + [-n, n]^d) \times \{0\} \textnormal{ to a } (y', t) \in (y + [-\tilde n, \tilde n]^d) \times \{ t \} )  >0 ,     \notag     
\end{align}
since the space-time box is finite.

By the definition of $K^\prime$ and by the domination of the fleas process and the dual process, there is an occupied cube of size $2n$ centered in $R_k$ and an occupied cube of size $ 2 \tilde n$ centered in $\tilde R_k$. By the duality, the probability that the two processes intersect at time $T/2$ is bounded from below by the LHS of \eqref{intersect inf}. Moreover, as long as $|k- l| > 1$, those intersection events between $R_k$ and $\tilde R_k$ and between $R_l$ and $\tilde R_l$ are independent. 

Let $N_2$ be such that at least one among $N_2$ independent trials whose success probability is equal to the LHS of \eqref{intersect inf} is successful with probability $\geq 1 - \varepsilon$. Thus, if $T$ is large enough such that the size of $K^\prime$ is $ \geq 2N_2$ with probability at least $1 - 3 \varepsilon$, then $B^B_\frac{T}{2}\cap B^ {T, D}     _\frac{T}{2} \neq\emptyset$ with high probability.


Now, we see that it is necessary that one of the following bad events must happen so that the event in \eqref{bad_event_bd} can happen. (1) $\xi_{n _  1} ^{(1)} = \emptyset$ or $\xi_{n _  2} ^{(2)} = \emptyset$ (with probability less than $2 \varepsilon$); (2) The cubes in step 2 are not present in time (with probability less than $3 \varepsilon$); (3) $|K| < 12N_2$ (with probability  less than $\varepsilon$ by the choice of $T$); (4) $|K^\prime| < |K| / 6$ (with probability $2\varepsilon$); (5) all independent trials as in \eqref{intersect inf} are unsuccessful (with probability $\varepsilon$). 

Summing up, the probability of bad events combined is less than $9\varepsilon$. By the arbitrary choice of $\varepsilon$, the proof is complete.
\end{proof}

With the previous lemma, we may now present the proof of Theorem \ref{cc}.
\begin{proof}[Proof of Theorem \ref{cc}]
We consider the following probability,
\begin{align}
 \nonumber
P (B_t^B \cap D\neq \emptyset) &=  P (B^B_{\frac t2}\cap B^{t, D } _\frac t2\neq\emptyset)\\
  &=P (B^B_{\frac t2}\neq \emptyset)P(B^{t, D }    _\frac t2 \neq \emptyset) - P(B^B _\frac{t}{2}\neq \emptyset, B^{t, D} _\frac{t}{2}\neq \emptyset, B^B_\frac{t}{2}\cap B^ {t, D}     _\frac{ t }{2}=\emptyset)\\
  \nonumber
  &=     P (B^B_{\frac t2}\neq \emptyset)P(B^{\frac t2, D }    _\frac t2 \neq \emptyset) - P(B^B _\frac{t}{2}\neq \emptyset, B^{t, D} _\frac{t}{2}\neq \emptyset, B^B_\frac{t}{2}\cap B^ {t, D}     _\frac{ t }{2}=\emptyset) .
\end{align}
As $t\to +\infty$, $P (B^B_{\frac t2}\neq \emptyset)     \to P (T_2 = \infty) $, $P(B^B _\frac{t}{2}\neq \emptyset, B^{t, D} _\frac{t}{2}\neq \emptyset, B^B_\frac{t}{2}\cap B^ {t, D}     _\frac{ t }{2}=\emptyset) \to 0$ by Lemma \ref{cc1}, and $P(B^{\frac t2, D }    _\frac t2 \neq \emptyset) =P (B^{\mathbb{Z}^d}_\frac t2 \cap D\neq \emptyset)  \to \bar \nu (\{E \subset \mathbb Z^d : E \cap D \neq     \emptyset \})$.  Thus,
\begin{equation}
   P (B_t^B \cap D\neq \emptyset) \to   P (T_2 = \infty)   \bar \nu (\{E \subset \mathbb Z^d : E \cap D \neq     \emptyset \}) 
   .
\end{equation}
Considering that the distribution of animals does not change, the proof is complete.
\end{proof}

\end {document}